\documentclass[11pt,twoside]{amsart}

\usepackage[all]{xy}
\usepackage{amsmath}
\usepackage{amsfonts}
\usepackage{amssymb}
\usepackage{amsthm}
\usepackage{enumerate}
\usepackage{color}
\usepackage{fullpage}
\usepackage{graphicx}
\usepackage{url}

\theoremstyle{definition}
\newtheorem{defn}{Definition}[section]

\newtheorem{question}[defn]{Question}
\theoremstyle{plain}
\newtheorem{thm}[defn]{Theorem}
\newtheorem{lem}[defn]{Lemma}
\newtheorem{prop}[defn]{Proposition}

\newtheorem{conj}[defn]{Conjecture}

\def\C{\ensuremath{\mathbb{C}}}

\def\P{\ensuremath{\mathbb{P}}}
\def\Q{\ensuremath{\mathbb{Q}}}
\def\R{\ensuremath{\mathbb{R}}}

\def\AA{\ensuremath{\mathcal A}}
\def\BB{\ensuremath{\mathcal B}}

\def\EE{\ensuremath{\mathcal E}}
\def\FF{\ensuremath{\mathcal F}}

\def\HH{\ensuremath{\mathcal H}}
\def\II{\ensuremath{\mathcal I}}

\def\NN{\ensuremath{\mathcal N}}
\def\OO{\ensuremath{\mathcal O}}

\def\TT{\ensuremath{\mathcal T}}

\def\ch{\mathop{\mathrm{ch}}\nolimits}

\def\Coh{\mathop{\mathrm{Coh}}\nolimits}

\def\Db{\mathop{\mathrm{D}^{\mathrm{b}}}\nolimits}

\def\End{\mathop{\mathrm{End}}\nolimits}
\def\lEnd{\mathop{\mathcal End}\nolimits}

\def\Hom{\mathop{\mathrm{Hom}}\nolimits}
\def\lHom{\mathop{\mathcal Hom}\nolimits}

\def\Tor{\mathop{\mathrm{Tor}}\nolimits}

\def\td{\mathop{\mathrm{td}}\nolimits}

\def\into{\ensuremath{\hookrightarrow}}
\def\onto{\ensuremath{\twoheadrightarrow}}

\begin{document}

\title{Bridgeland Stability on Blow Ups and Counterexamples}

%\author{Cristian Martinez}
%\address{University of California, Department of Mathematics, South Hall, Room 6607, 
%Santa Barbara, CA 93106-3080, USA}
%\email{martinez@math.ucsb.edu}
%\urladdr{http://web.math.ucsb.edu/~martinez/}

\author{Cristian Martinez}
\address{Universidad de Los Andes, Department of Mathematics, Cra 1 No 18A-12, 
Bogot\'a, Colombia}
\email{cm.martineze@uniandes.edu.co}
\urladdr{https://sites.google.com/site/cristianmathinez/}

\author{Benjamin Schmidt}
\address{The University of Texas at Austin, Department of Mathematics, 2515 Speedway, RLM 8.100, Austin, TX 78712, USA}
\email{schmidt@math.utexas.edu}
\urladdr{https://sites.google.com/site/benjaminschmidtmath/}

%Omprokash Das address

\address{Department of Mathematics\\
University of California, Los Angeles\\
520 Portola Plaza\\
Math Sciences Building 6363.}
\email{das@math.ucla.edu}
\urladdr{https://www.math.ucla.edu/~das/}



\begin{abstract}
We give further counterexamples to the conjectural construction of Bridgeland stability on threefolds due to Bayer, Macr\`i, and Toda. This includes smooth projective threefolds containing a divisor that contracts to a point, and Weierstra{\ss} elliptic Calabi-Yau threefolds. Furthermore, we show that if the original conjecture, or a minor modification of it, holds on a smooth projective threefold, then the space of stability conditions is non-empty on the blow up at an arbitrary point. More precisely, there are stability conditions on the blow up for which all skyscraper sheaves are semistable.
\end{abstract}

\maketitle

\setcounter{tocdepth}{1}
\tableofcontents

\section{Introduction}

Since Bridgeland introduced stability conditions on triangulated categories in \cite{Bri07:stability_conditions}, the topic has been haunted. He was motivated by the attempt to further understand homological mirror symmetry related to Calabi-Yau threefolds, but to this day we do not know how to construct stability conditions on threefolds in general. While the theory has flourished with many applications on curves, surfaces, and quiver representations, the threefold problem has persisted.

In \cite{BMT14:stability_threefolds} Bayer, Macr\`i, and Toda proposed a conjectural construction of Bridgeland stability on threefolds (see Section \ref{sec:prelim} for full details). They define the intermediate notion of tilt stability analogously to the construction of Bridgeland stability in the surface case. Let $X$ be a smooth projective threefold with ample polarization $H$, and let $B_0$ be an arbitrary $\Q$-divisor on $X$. For each $\beta \in \R$ one tilts the category of coherent sheaves to obtain the new heart of a bounded t-structure $\Coh^{\beta}(X)$ consisting of certain two-term complexes. After fixing another real number $\alpha > 0$, the new slope function is given by
\[
\nu_{\alpha, \beta} := \frac{H \cdot \ch^B_2 - \frac{\alpha^2}{2} H^3 \cdot \ch^B_0}{H^2 \cdot \ch^B_1},
\]
where $\ch^{B} = e^{-B} \cdot \ch$, and $B = B_0 + \beta H$. Note that if $B$ is an integral divisor, then we have $\ch^{B}(E) = \ch(E \otimes \OO(-B))$. In order to construct Bridgeland stability, they propose another tilt of the category $\Coh^{\beta}(X)$. They managed to prove all necessary properties except for the following conjecture.

\theoremstyle{plain}
\newtheorem*{conj:bmt}{Conjecture \ref{conj:bmt}}
\begin{conj:bmt}[{\cite[Conjecture 1.3.1]{BMT14:stability_threefolds}}]
For any $\nu_{\alpha,\beta}$-stable object $E \in \Coh^{\beta}(X)$ with $\nu_{\alpha,\beta}(E)=0$ the inequality
\[
\ch^B_3(E) \leq \frac{\alpha^2}{6} H^2 \cdot \ch^B_1(E)
\]
holds.
\end{conj:bmt}

%\footnote{\textcolor{red}{These paragraphs appear in section 2 as well, we should sumarize the information either here or there}}
This inequality was first proved for $\P^3$ in \cite{Mac14:conjecture_p3}. It was shown to hold for the smooth quadric in $\P^4$ in \cite{Sch14:conjecture_quadric} and was later generalized to all Fano threefolds of Picard rank one in \cite{Li15:conjecture_fano_threefold}. The case of abelian threefolds was settled with two independent proofs \cite{BMS14:abelian_threefolds,MP16:conjecture_abelian_threefoldsII}. Most recently, it was shown in \cite{Kos17:stability_products} for the case of $\P^2 \times E$, $\P^1 \times \P^1 \times E$, and $\P^1 \times A$, where $E$ is an arbitrary elliptic curve and $A$ is an arbitrary abelian surface.

It turns out that the conjecture does not hold in general which was first pointed out for the blow up of $\P^3$ at a point in \cite{Sch17:counterexample}. Moreover, \cite{Kos17:stability_products} gave a counterexample for Calabi Yau threefolds containing a plane. A modified conjecture was proved for all Fano threefolds in \cite{BMSZ16:stability_fano} and \cite{Piy17:Fano}. This answers the following question affirmatively in that case.

\theoremstyle{definition}
\newtheorem*{q:new_conjecture}{Question \ref{q:new_conjecture}}
\begin{q:new_conjecture}[\cite{BMSZ16:stability_fano}]
Is there a cycle $\Gamma\in A_1(X)_\Q$ depending at most on $H$ and $B_0$ such that $\Gamma \cdot H \geq 0$ and for any $\nu_{\alpha,\beta}$-stable object $E$ with $\nu_{\alpha,\beta}(E)=0$, we have 
\[
\ch_3^B(E) \leq \Gamma \cdot \ch_1^B(E) + \frac{\alpha^2}{6} H^2 \cdot \ch_1^B(E)?
\]
\end{q:new_conjecture}

Correcting and proving the conjecture is fundamental to the advancement of the theory of Bridgeland stability on threefolds. We believe that the following additional counterexamples should shed more light on this question.

\begin{thm}
Let $f: X \to Y$ be a birational morphism between projective threefolds, where $X$ is smooth. Let $D \subset X$ be an effective divisor such that $-D$ is $f$-ample and $f(D)$ is a point. Then there is a polarization $H$ on $X$ such that Conjecture \ref{conj:bmt} fails for $\OO_D$ with $B_0 = 0$ for some $(\alpha, \beta)$.

In particular, Conjecture \ref{conj:bmt} fails for any Weierstra{\ss} elliptic Calabi-Yau threefold over a del Pezzo surface.
\end{thm}

The condition on $-D$ being $f$-ample is very weak, and such a divisor exists for any birational divisorial contraction, where $Y$ is normal and $\Q$-factorial. In order to prove this theorem, we first establish the technical Lemma \ref{counterex:criterion} that allows to determine whether the conjecture fails on the structure sheaf of a divisor. In case of the Weierstra{\ss} elliptic Calabi-Yau threefold, we get the counterexample from the structure sheaf of the image of its canonical section. Omprokash Das pointed out to us that in the case of the Weierstra{\ss} elliptic Calabi-Yau threefold there is a birational morphism contracting exactly the image of the canonical section to a point, and therefore, it is a special case of the theorem. This is proven in Appendix A.

We would like to point out that the counterexample for Calabi-Yau threefolds containing a plane in \cite{Kos17:stability_products} is also a consequence of Lemma \ref{counterex:criterion}. Indeed, as pointed out by Koseki, if a Calabi-Yau threefold contains a plane, then such a plane can be contracted and 
its structure sheaf contradicts Conjecture \ref{conj:bmt} just as in Section \ref{contracting:divisors}. All together this suggests an unexpected relation between the construction of Bridgeland stability on smooth projective threefolds and their birational geometry.

It is noteworthy that there is still no known counterexample in Picard rank one. This allows to draw a comparison to the classical Bogomolov-Gieseker inequality on surfaces for semistable sheaves. It says that any torsion-free semistable sheaf $E$ satisfies the inequality
\[
\ch_1(E)^2 - 2 \ch_0(E) \ch_2(E) \geq 0.
\]
This inequality might fail for structure sheaves of divisors, too. However, as in our case this does not happen in the case of Picard rank one.

Finally, we further investigate the blow up case in Section \ref{sec:blow_up}.

\begin{thm}
Assume that $X$ is a smooth projective threefold where Question \ref{q:new_conjecture} has an affirmative answer for some polarization $H$ and $\Q$-divisor $B_0$. Then the induced upper half-plane of stability conditions on $X$ embeds into the space of stability conditions on the blow up of $X$ at an arbitrary point.

For these stability conditions skyscraper sheaves $\C(x)$ on the blow up are all semistable. A skyscraper sheaf $\C(x)$ is stable if and only if $x$ does not lie on the exceptional divisor.
\end{thm}

Let $f: \tilde{X} \to X$ be the blow up. As polarization we choose the pullback $\tilde{H} = f^*H$ which is not ample anymore, but just nef. This makes a modification of the construction of the hearts of bounded t-structures necessary (see Section \ref{sec:blow_up} for details). Within this slightly modified framework, the class $\tilde{\Gamma}$ on the blow up is given by $f^* \Gamma - \tfrac{1}{6}E^2$, and the $\Q$-divisor $\tilde{B}_0$ is $f^*B_0 + 2E$.

The proof of this Theorem is based on a result by Toda from \cite{Tod13:extremal_contractions}. He proved that the derived category $\Db(\tilde{X})$ is equivalent to the bounded derived category of the abelian category of finitely generated $\BB$-modules $\Db(X,\BB)$ for a certain sheaf of finitely generated $\OO_X$-algebras $\BB$. We carefully study preservation of stability for the forgetful functor $\Db(X, \BB) \to \Db(X)$ to construct Bridgeland stability on $\Db(X, \BB) \cong \Db(\tilde{X})$ (see Lemma \ref{lem:mu_preserved_forgetful} and Lemma \ref{lem:nu_preserved_forgetful}).

\subsection*{Acknowledgments}
We would like to thank Arend Bayer, Aaron Bertram, Omprokash Das, Jason Lo, Emanuele Macr\`i, and David R. Morrison for discussions on the topic of this article. We also thank the referees for carefully reading the article.

\subsection*{Notation}
\begin{center}
  \begin{tabular}{ r l }
    $X$ & smooth projective threefold over $\C$ \\
    $H$ & fixed ample divisor on $X$ \\
    $\Db(X)$ & bounded derived category of coherent \\ & sheaves on $X$ \\
    $\ch(E)$ & Chern character of an object $E \in \Db(X)$  \\
    $\ch_{\leq l}(E)$ & $(\ch_0(E), \ldots, \ch_l(E))$ \\
    $H \cdot \ch(E)$ & $(H^3 \cdot \ch_0(E), H^2 \cdot \ch_1(E), H \cdot \ch_2(E) 
    \ch_3(E))$ \\
    $H \cdot \ch_{\leq l}(E)$ & $(H^3 \cdot
    \ch_0(E), \ldots, H^{3-l} \cdot \ch_l(E))$
  \end{tabular}
\end{center}

%%%%%%%%%%%%%%%%%%%%%%%%%%%%%%%%%%%%%%%%%%%%%%%%%
%%%%%%%%%%%%%%%%%%%%%%%%%%%%%%%%%%%%%%%%%%%%%%%%%
%%%%%%%%%%%%%%%%%%%%%%%%%%%%%%%%%%%%%%%%%%%%%%%%%
\section{Preliminaries}
\label{sec:prelim}

Throughout this section, we fix a smooth projective threefold $X$, an integral ample divisor class $H$, and an arbitrary $\Q$-divisor class $B_0$. Moreover, for real numbers $\alpha > 0$ and $\beta \in \R$, we define $\omega = \alpha H$ and $B = B_0 + \beta H$. While $X$, $H$, and $B_0$ are fixed, we view $\alpha$ and $\beta$ as varying parameters. The goal is to construct an upper half-plane of stability conditions based on these parameters.

The \emph{classical slope} for a coherent sheaf $E \in \Coh(X)$ is defined as
\[
\mu(E) := \frac{H^2\cdot \ch_1(E)}{H^3 \cdot \ch_0(E)},
\]
where as usual division by zero is interpreted as $+\infty$. A coherent sheaf $E$ is called \emph{slope (semi)stable} if for any non trivial proper subsheaf $F \subset E$ the inequality $\mu(F) < (\leq) \mu(E/F)$ holds.

To ease notation, we define the \emph{twisted Chern character} $\ch^B$ as $e^{-B} \cdot \ch$. Note that in the case where $B$ is integral, we simply have $\ch^B(E) = \ch(E \otimes \OO(-B))$ This definition expands to
\begin{align*}
\ch^B_0 &= \ch_0, \\
\ch^B_1 &= \ch_1 - B \cdot \ch_0 ,\\
\ch^B_2 &= \ch_2 - B \cdot \ch_1 + \frac{B^2}{2} \cdot \ch_0, \\
\ch^B_3 &= \ch_3 - B \cdot \ch_2 + \frac{B^2}{2} \cdot \ch_1 - \frac{B^3}{6} \cdot \ch_0.
\end{align*}

In the case $B_0 = 0$, we write $\ch^{\beta} := \ch^{B}$. The theory of tilting is used to construct another heart of a bounded t-structure. For more information on the general method of tilting we refer to \cite{HRS96:tilting} and \cite{BvdB03:functors}. A \emph{torsion pair} on the category of coherent sheaves can be defined by
\begin{align*}
\TT_{\beta} &:= \left\{E \in \Coh(X) : \text{any quotient $E \onto G$ satisfies $\mu(G) > \frac{H^2 \cdot B}{H^3}$} \right\}, \\
\FF_{\beta} &:= \left\{E \in \Coh(X) : \text{any non-trivial subsheaf $F \subset E$ satisfies $\mu(F) \leq \frac{H^2 \cdot B}{H^3}$} \right\}.
\end{align*}
A new heart of a bounded t-structure is given as the extension closure $\Coh^{\beta}(X) := \langle \FF_{\beta}[1],\TT_{\beta} \rangle$. This means objects in $\Coh^{\beta}(X)$ are given by morphisms between coherent sheaves with kernel in $\FF_{\beta}$ and cokernel in $\TT_{\beta}$. The \emph{tilt slope} is defined as
\[
\nu_{\alpha, \beta} := \frac{H \cdot \ch^B_2 - \frac{\alpha^2}{2} H^3 \cdot \ch^B_0}{H^2 \cdot \ch^B_1}.
\]
As before, an object $E \in \Coh^{\beta}(X)$ is called \emph{tilt-(semi)stable} (or \emph{$\nu_{\alpha,\beta}$-(semi)stable}) if for any non trivial proper subobject $F \subset E$ the inequality $\nu_{\alpha, \beta}(F) < (\leq) \nu_{\alpha, \beta}(E/F)$ holds.

\begin{thm}[{Bogomolov Inequality for Tilt Stability, \cite[Corollary 7.3.2]{BMT14:stability_threefolds}}]
\label{thm:bogomolov}
Any $\nu_{\alpha, \beta}$-semistable object $E \in \Coh^{\beta}(X)$ satisfies
\begin{align*}
\Delta(E) &= (H^2 \cdot \ch_1^B(E))^2 - 2(H^3 \cdot \ch_0^B(E))(H \cdot \ch_2^B(E)) \\
&= (H^2 \cdot \ch_1^{B_0}(E))^2 - 2(H^3 \cdot \ch^{B_0}_0(E))(H \cdot \ch^{B_0}_2(E)) \geq 0.
\end{align*}
\end{thm}

Let $\Lambda$ be the smallest lattice containing the image of the map $H \cdot \ch^{B_0}_{\leq 2}$. Clearly, $\Lambda$ has rank three. For any $i \in \{0, 1, 2\}$ there is a unique function $( \cdot )_i: \Lambda \to \Q$ such that for $v = H \cdot \ch^{B_0}_{\leq 2}(E)$, we have $v_i = H^{3 - i} \cdot \ch^{B_0}_i(E)$.
%The image of the map $H \cdot \ch^{B_0}_{\leq 2}$ is contained in $\Lambda = \Z \oplus \Z \oplus \tfrac{1}{2} \Z$.
Notice that $\nu_{\alpha, \beta}$ factors through $H \cdot \ch^{B_0}_{\leq 2}$. Varying $(\alpha, \beta)$ changes the set of (semi)stable objects. A \emph{numerical wall} in tilt stability with respect to a class $v \in \Lambda$ is a non trivial proper subset $W$ of the upper half plane given by an equation of the form $\nu_{\alpha, \beta}(v) = \nu_{\alpha, \beta}(w)$ for another class $w \in \Lambda$. A subset $S$ of a numerical wall $W$ is called an \emph{actual wall} if the set of semistable objects with class $v$ changes at $S$. The structure of walls in tilt stability is rather simple. This is sometimes also called Bertram's Nested Wall Theorem and the first full proof appears in \cite{Mac14:nested_wall_theorem}
%Part (1) - (5) is usually called Bertram's Nested Wall Theorem and appeared first in \cite{Mac14:nested_wall_theorem}, while part (6), (7), and (8) are in Lemma 2.7 and Appendix A of \cite{BMS14:abelian_threefolds}. The last part of (8) about reflexivity is to be found in \cite[Proposition 3.1]{LM16:examples_tilt}.

\begin{thm}[Structure Theorem for Tilt Stability]
\label{thm:Bertram}
Let $v \in \Lambda$ be a fixed class. All numerical walls in the following statements are with respect to $v$.
\begin{enumerate}
  \item Numerical walls in tilt stability are either semicircles with center on the $\beta$-axis or rays parallel to the $\alpha$-axis.
  \item If two numerical walls given by classes $w,u \in \Lambda$ intersect, then $v$, $w$ and $u$ are linearly dependent. In particular, the two walls are completely identical.
  \item The curve $\nu_{\alpha, \beta}(v) = 0$ is given by a hyperbola, or a vertical line in the special case $v_0 = 0$. Moreover, it intersects all semicircular walls at their top point.
  \item If $v_0 \neq 0$, there is exactly one numerical vertical wall given by $\beta = v_1/v_0$. If $v_0 = 0$, there is no actual vertical wall.
  \item If a numerical wall has a single point at which it is an actual wall,
  then all of it is an actual wall.
%  \item If there is an actual wall numerically defined by an exact sequence of tilt semistable objects $0 \to F \to E \to G \to 0$ such that $H \cdot \ch_{\leq 2}(E) = v$, then 
%  \[
%  \Delta(F) + \Delta(G) \leq \Delta(E).
%  \]
%  Moreover, equality holds if and only if either $H \cdot \ch_{\leq 2}(F) = 0$, $H \cdot \ch_{\leq 2}(G) = 0$, or both $\Delta(E) = 0$ and $H \cdot \ch_{\leq 2}(F)$, $H \cdot \ch_{\leq 2}(G)$, and $H \cdot \ch_{\leq 2}(E)$ are all proportional.
%  \item If $\Delta(E) = 0$ for a tilt semistable object $E$, then $E$ can only be destabilized at the unique numerical vertical wall.
%  \item If $E$ is a tilt stable object for fixed $\beta \in \R$ and all $\alpha \gg 0$, then $E$ is one of the following.
%  \begin{itemize}
%  \item If $\ch_0(E) \geq 0$, then $E$ is a slope semistable sheaf.
%  \item If $\ch_0(E) < 0$, then $H^0(E)$ is a torsion sheaf supported in dimension smaller than or equal to $1$ and $H^{-1}(E)$ is a reflexive slope semistable sheaf with positive rank.
% \end{itemize}
\end{enumerate}
\end{thm}

The following lemma is useful in computations and first appeared in \cite{CH15:nef_cones}. We refer to \cite[Lemma 7.2]{MS16:lectures_notes} for a simple proof.

\begin{lem}
\label{lem:higherRankBound}
Let $0 \to F \to E \to G \to 0$ be an exact sequence in $\Coh^{\beta}(X)$ defining a non empty semicircular wall $W$. Assume further that $\ch_0(F) > \ch_0(E) \geq 0$. Then the radius $\rho_W$ satisfies the inequality
\[
\rho_W^2 \leq \frac{\Delta(E)}{4 H^3 \cdot \ch_0(F) (H^3 \cdot \ch_0(F) - H^3 \cdot \ch_0(E))}.
\]
\end{lem}

%Another simple but yet important observation is that since $\nu_{\alpha, \beta}$ does not involve the third Chern character, then the subcategory of 0-dimensional sheaves remains unaffected after tilting. More precisely, we have the following lemma.
Note that by definition an object in $\Coh^{\beta}(X)$ is supported in dimension zero if and only if it is a zero dimensional torsion sheaf. A simple but yet important observation is that their subcategory in $\Coh^{\beta}(X)$ is closed under taking subobjects and quotients. More precisely, we have the following lemma.

\begin{lem}
\label{tor:quotient}
Let $0\to K \to A \to B \to 0$ be a short exact sequence in $\Coh^{\beta}(X)$ and suppose that $A$ is a zero-dimensional sheaf, then so are $K$ and $B$.
\end{lem}

\begin{proof}
Consider the long exact sequence of cohomologies in $\Coh(X)$
\[
0 \to \HH^{-1}(B) \to K \to A \to \HH^{0}(B) \to 0.
\]
Thus, $\HH^{0}(B)$ is a zero-dimensional sheaf as well. If $\HH^{-1}(B)$ is nonzero, then we would have $\mu(\HH^{-1}(B)) = \mu(K)$ contradicting the definition of $\Coh^{\beta}(X)$. 
\end{proof}

A generalized Bogomolov type inequality involving third Chern characters for tilt semistable objects with $\nu_{\alpha, \beta} = 0$ has been conjectured in \cite{BMT14:stability_threefolds}. Its main goal was the construction of Bridgeland stability conditions on arbitrary threefolds.

\begin{conj}[{\cite[Conjecture 1.3.1]{BMT14:stability_threefolds}}]
\label{conj:bmt}
For any $\nu_{\alpha,\beta}$-stable object $E \in \Coh^{\beta}(X)$ with $\nu_{\alpha,\beta}(E)=0$ the inequality
\[
\ch^B_3(E) \leq \frac{\alpha^2}{6} H^2 \cdot \ch^B_1(E)
\]
holds.
\end{conj}

It turns out that the conjecture does not hold in general which was first pointed out for the blow up of $\P^3$ at a point in \cite{Sch17:counterexample}. We give further counterexamples in the next section. A modified conjecture was proved for all Fano threefolds in \cite{BMSZ16:stability_fano} and \cite{Piy17:Fano}. This answers the following question affirmatively in this case.

% \footnote{\textcolor{red}{These same two paragraphs appear in the introduction, should we rewrite them in a different way? Maybe summarize both paragraphs in one?}}This inequality was first proved for $\P^3$ in \cite{Mac14:conjecture_p3}. It was shown to hold for the smooth quadric in $\P^4$ in \cite{Sch14:conjecture_quadric}, and later generalized to all Fano threefolds of Picard rank one in \cite{Li15:conjecture_fano_threefold}. The case of abelian threefolds was settled with two independent proofs \cite{BMS14:abelian_threefolds,MP16:conjecture_abelian_threefoldsII}. Most recently, it was shown in \cite{Kos17:stability_products} for the case of $\P^2 \times E$, $\P^1 \times \P^1 \times E$, and $\P^1 \times A$, where $E$ is an arbitrary elliptic curve and $A$ is an arbitrary abelian surface.

% It turns out that the conjecture does not hold in general which was first pointed out for the blow up of $\P^3$ at a point in \cite{Sch17:counterexample}. Moreover, \cite{Kos17:stability_products} gave a counterexample for Calabi Yau threefolds containing a plane.

\begin{question}[\cite{BMSZ16:stability_fano}]
\label{q:new_conjecture}
Is there a cycle $\Gamma\in A_1(X)_\R$ depending at most on $H$ and $B_0$ such that $\Gamma \cdot H \geq 0$ and for any $\nu_{\alpha,\beta}$-stable object $E$ with $\nu_{\alpha,\beta}(E)=0$, we have 
\[
\ch_3^B(E) \leq \Gamma \cdot \ch_1^B(E) + \frac{\alpha^2}{6} H^2 \cdot \ch_1^B(E)?
\]
\end{question}

The same way as in \cite{BMT14:stability_threefolds} one can use an affirmative answer to this question to construct Bridgeland stability as follows. One repeats the process of tilting by replacing $\Coh(X)$ with $\Coh^{\beta}(X)$ and $\mu$ with $\nu_{\alpha, \beta}$. A torsion pair on $\Coh^{\beta}(X)$ is then defined by
\begin{align*}
\TT'_{\alpha, \beta} &:= \left\{E \in \Coh^{\beta}(X) : \text{any quotient $E \onto G$ satisfies $\nu_{\alpha, \beta}(G) > 0$} \right\}, \\
\FF'_{\alpha, \beta} &:= \left\{E \in \Coh^{\beta}(X) : \text{any non-trivial subobject $F \subset E$ satisfies $\nu_{\alpha, \beta}(F) \leq 0$} \right\}.
\end{align*}
The heart of a bounded t-structure is given by the extension closure $\AA^{\alpha, \beta}(X) := \langle \FF'_{\alpha, \beta}[1], \TT'_{\alpha, c\beta} \rangle$. As is customary with Bridgeland stability, one defines the following \emph{central charge} instead of a slope that depends on an additional parameter $s > \tfrac{1}{6}$
\[
Z^{\Gamma}_{\alpha, \beta, s} = -\ch_3^B + s \alpha^2 H^2 \cdot \ch_1^B + \Gamma \cdot \ch_1^B + i \left( H \cdot \ch_2^B - \frac{\alpha^2}{2} H^3 \cdot \ch_0^B \right).
\]
The corresponding \emph{Bridgeland slope} is then given by
\[
\lambda^{\Gamma}_{\alpha, \beta,s} := \frac{\ch_3^B - s \alpha^2 H^2 \cdot \ch_1^B - \Gamma \cdot \ch_1^B}{H \cdot \ch_2^B - \frac{\alpha^2}{2} H^3 \cdot \ch_0^B}.
\]

With the same arguments as in \cite{BMT14:stability_threefolds} the pair $(\AA^{\alpha, \beta}, Z^{\Gamma}_{\alpha, \beta, s})$ is a Bridgeland stability condition on $\Db(X)$ due to the inequality in Question \ref{q:new_conjecture}. This inequality implies the most critical property of a Bridgeland stability condition, namely for any $E \in \AA^{\alpha, \beta}(X)$ with $\Im Z^{\Gamma}_{\alpha, \beta, s}(E) = 0$ one gets $\Re Z^{\Gamma}_{\alpha, \beta, s}(E) < 0$.

%%%%%%%%%%%%%%%%%%%%%%%%%%%%%%%%%%%%%%%%%%%%%%%%%
%%%%%%%%%%%%%%%%%%%%%%%%%%%%%%%%%%%%%%%%%%%%%%%%%
%%%%%%%%%%%%%%%%%%%%%%%%%%%%%%%%%%%%%%%%%%%%%%%%%
\section{Counterexamples}
As shown by the first author in \cite{Sch17:counterexample} the generalized Bogomolov-Gieseker inequality (Conjecture \ref{conj:bmt}) fails on the blow-up of $\mathbb{P}^3$ at one point. The purpose of this section is to explore a more general class of counterexamples in hopes of shedding some light on how to choose the class $\Gamma$ of Question \ref{q:new_conjecture}.

Our approach is very simple. In the same way as the structure sheaf of a curve of negative self intersection violates the usual Bogomolov-Gieseker inequality on a surface, we want to look for divisors on a threefold with special intersection properties so that their structure sheaves will violate Conjecture \ref{conj:bmt}.

As before, let $X$ be a smooth projective threefold, $H$ an integral ample class, and $B=\beta H$ an $\mathbb{R}$-divisor on $X$. For an effective integral divisor $D$ on $X$ we can compute the twisted Chern characters
\begin{align*}
\ch^{\beta}_0(\OO_D) &= 0,\\
\ch^{\beta}_1(\OO_D) &= D,\\
\ch^{\beta}_2(\OO_D) &= -\frac{D^2}{2} - \beta H\cdot D,\\
\ch^{\beta}_3(\OO_D) &= \frac{D^3}{6} + \beta H \cdot \frac{D^2}{2} + \frac{\beta^2}{2} H^2 \cdot D.
\end{align*} 
To produce a counterexample to Conjecture \ref{conj:bmt} we want to find $\beta$ and $\alpha$ so that:
\begin{enumerate}[(a)]
\item $\displaystyle \OO_D$ is $\displaystyle \nu_{\alpha,\beta}$-semistable,
\item $\displaystyle \nu_{\alpha,\beta}(\OO_D)=0$, and
\item $\displaystyle \ch^{\beta}_3(\OO_D) - \frac{\alpha^2}{6} H^2 \cdot \ch^{\beta}_1(\OO_D) > 0$.
\end{enumerate}
Notice that there is only one value $\beta_0$ for which condition (a) is satisfied. Indeed, $\nu_{\alpha,\beta}(\OO_D) = 0$ if and only if
\[
\beta = \beta_0 = -\frac{D^2 \cdot H}{2D \cdot H^2}.
\] 
We have
\begin{align*}
\ch^{\beta_0}_3(\OO_D) - \frac{\alpha^2}{6} H^2 \cdot \ch^{\beta_0}_1(\OO_D) &= \frac{D^3}{6} + \beta_0 \frac{D^2 \cdot H}{2} + \beta_0^2 \frac{D \cdot H^2}{2} - \frac{\alpha^2}{6} D \cdot H^2\\
&= \frac{D^3}{6} - \frac{(D^2 \cdot H)^2}{4D \cdot H^2} + \frac{(D^2 \cdot H)^2 (D \cdot H^2)}{8(D \cdot H^2)^2} - \frac{\alpha^2}{6} D \cdot H^2\\
&= \frac{D^3}{6} - \frac{1}{8} \frac{(D^2 \cdot H)^2}{D \cdot H^2} - \frac{\alpha^2}{6} D \cdot H^2.
\end{align*}
Thus, condition (c) is satisfied if and only if
\begin{equation}
\label{condition(c)}
\frac{\alpha^2}{6} D \cdot H^2 < \frac{D^3}{6} - \frac{1}{8} \frac{(D^2 \cdot H)^2}{D \cdot H^2}.
\end{equation}
We are left to find a range of values for $\alpha$ so that condition (a) is satisfied. First of all, notice that since $\OO_D$ is a Gieseker semistable torsion sheaf, $\OO_D$ is $\nu_{\alpha,\beta}$-semistable for all $\beta$, and $\alpha\gg 0$. Moreover, the walls for tilt semistability in the $(\alpha,\beta)$-plane for the Chern character $\ch(\OO_D)$ are semicircles with center $(0,\beta_0)$. 

By Theorem \ref{thm:Bertram}, we know that if $\OO_D$ is destabilized we must have a short exact sequence 
$0 \to A \to \OO_D \to B \to 0$ along a semicircular wall $W$ of radius $\rho_W$. Note that the point $(\alpha, \beta_0) \in W$ is given by $\alpha = \rho_W$. Since any subobject of a sheaf in $\Coh^{\beta}(X)$ is a sheaf, the rank of $A$ is non-negative. If $A$ had rank zero, than a straightforward computation shows that $A$ either destabilizes $\OO_D$ for all $(\alpha, \beta)$ or none, a contradiction. Thus, $A$ must have rank at least one, and we can use Lemma \ref{lem:higherRankBound} to get
\[
\rho_W^2 \leq \frac{\Delta(\OO_D)}{4 (H^3 \cdot \ch_0(F))^2} \leq \frac{(D\cdot H^2)^2}{4(H^3)^2}.
\]

In particular, this shows that $\OO_D$ is $\nu_{\alpha,\beta_0}$-semistable for 
\begin{equation}
\label{condition(a)}
\alpha \geq \frac{D\cdot H^2}{2H^3}.
\end{equation}
Combining \eqref{condition(c)} and \eqref{condition(a)} we obtain the following Lemma.

\begin{lem}
\label{counterex:criterion}
Let $X$ be a smooth projective threefold. Suppose that there is an effective divisor $D$ and an ample divisor $H$ such that
\begin{equation}
\label{counterexample}
D^3 > \frac{(D \cdot H^2)^3}{4(H^3)^2} + \frac{3}{4} \frac{(D^2 \cdot H)^2}{D \cdot H^2}.
\end{equation}
Then there exists a pair $(\alpha_0, \beta_0)$ such that $\OO_D$ violates Conjecture \ref{conj:bmt}.
\end{lem}

\subsection{Contracting divisors}\label{contracting:divisors}

Let $X$ be a smooth projective threefold and suppose that there is a projective morphism $\pi \colon X \to Y$ with exceptional locus a divisor $D$ that gets contracted to a point by $\pi$ and such that $-D$ is relatively ample. Let $A$ be an ample divisor on $Y$, and let $L = \pi^*A$. Then  for $m \gg 0$ the divisor $H = mL - D$ is ample on $X$. Since $L \cdot D = 0$ then one can easily compute 
\[
H^3 = m^3 L^3 - D^3,\ D \cdot H^2 = D^3,\ \text{and}\  D^2 \cdot H =-D^3.
\]
Notice that
\begin{align*}
D^3 - \frac{(D \cdot H^2)^3}{4(H^3)^2} - \frac{3}{4} \frac{(D^2 \cdot H)^2}{D \cdot H^2} &= D^3-\frac{(D^3)^3}{4(m^3 L^3 - D^3)^2} - \frac{3}{4} \frac{(-D^3)^2}{D^3}\\
%&= D^3-\frac{(D^3)^3}{4(m^3L^3-D^3)^2}-\frac{3}{4}D^3\\
&= \frac{D^3}{4} - \frac{(D^3)^3}{4(m^3 L^3 - D^3)^2}.
\end{align*}
This quantity is positive for $m \gg 0$ since $-D$ is relatively ample and so $D^3 > 0$. Therefore by Lemma \ref{counterex:criterion}, $\OO_D$ violates Conjecture \ref{conj:bmt}. In particular, Conjecture \ref{conj:bmt} fails on any blow up of a smooth threefold at a point.

\subsection{Weierstra{\ss} threefolds}

Let $p \colon X \to S$ be a smooth elliptic Calabi-Yau threefold over a del Pezzo surface $S$ in Weierstra{\ss} form with canonical section $\sigma \colon S \to X$. We refer to \cite[Section 6.2]{BBH09:fm_nahm_transform} for additional background on these threefolds.

Let us denote by $\Theta \subset X$ the image of $\sigma$ and by $K_S$ the canonical divisor of $S$. Then as observed by Diaconescu in \cite[Corollary 2.2(i)]{Dia16:fm_elliptic_CY3}, a divisor class
\[
H=t\Theta + p^*\eta
\]
with $t>0$ is ample on $X$ if and only if $\eta + tK_S$ is ample on $S$. For any $t >0$, we can set $\eta = -(1+t) K_S$. Using that $S$ is del Pezzo, we get that the divisor $H = t\Theta - (1+t) p^*K_S$ is ample. Using the adjunction formula for $\Theta \hookrightarrow X$ we obtain 
\[
\Theta^2 = \Theta \cdot p^*(K_S).
\]
Then one can compute the intersection numbers
\[
\Theta^3 = K_S^2 =\Theta \cdot H^2,\ \ \Theta^2 \cdot H = -K_S^2,\ \  \text{and}\ H^3= (t^3 + 3t^2 + 3t) K_S^2.
\]
Thus,
\[
\Theta^3 - \frac{(\Theta \cdot H^2)^3}{4(H^3)^2} - \frac{3}{4} \frac{(\Theta^2 \cdot H)^2}{\Theta \cdot H^2} = \frac{K_S^2}{4} \left(1 - \frac{1}{(t^3 + 3t^2 + 3t)^2} \right),\\
%&>0\ \ \text{for}\ \ t>2^{1/3}-1.
\]
which is positive for $t > 2^{1/3}-1$. Therefore, Lemma \ref{counterex:criterion} implies that Conjecture \ref{conj:bmt} fails for all smooth Weierstra{\ss} Calabi-Yau threefolds over a del Pezzo surface.

%%%%%%%%%%%%%%%%%%%%%%%%%%%%%%%%%%%%%%%%%%%%%%%%%
%%%%%%%%%%%%%%%%%%%%%%%%%%%%%%%%%%%%%%%%%%%%%%%%%
%%%%%%%%%%%%%%%%%%%%%%%%%%%%%%%%%%%%%%%%%%%%%%%%%
\section{Blowing up a point}
\label{sec:blow_up}

Let $X$ be a smooth projective threefold, $P$ be a point on $X$, $f: \tilde{X} \to X$ be the blow up of $X$ at $P$, and $E$ be the exceptional divisor. We will construct Bridgeland stability conditions on $\tilde{X}$ provided that the generalized Bogomolov-Gieseker inequality Conjecture \ref{conj:bmt} holds on $X$ or more generally Question \ref{q:new_conjecture} has an affirmative answer answer on $X$.

Let $H$ be an ample divisor on $X$ and $B_0$ any $\Q$-divisor on $X$. By assumption there is a cycle $\Gamma \in A_1(X)_{\Q}$ depending at most on $H$ and $B_0$ such that $\Gamma \cdot H \geq 0$ and for any $\nu_{\alpha,\beta}$-stable object $E$ with $\nu_{\alpha,\beta}(E)=0$, we have 
\[
\ch_3^B(E) \leq \Gamma \cdot \ch_1^B(E) + \frac{\alpha^2}{6} H^2 \cdot \ch_1^B(E).
\]
We define corresponding classes on $\tilde{X}$ as follows
\begin{align*}
\tilde{H} &= f^*H, \\
\tilde{B}_0 &= f^*B_0 + 2E, \\
\tilde{\Gamma} &= f^* \Gamma - \frac{E^2}{6}, \\
\tilde{\omega} &= \alpha \tilde{H}, \\
\tilde{B} &= \tilde{B_0} + \beta \tilde{H},
\end{align*}
where as previously $\alpha, \beta \in \R$ with $\alpha > 0$. For any $s > \tfrac{1}{6}$, we define a function on $\Db(\tilde{X})$ as
\[
Z_{\alpha, \beta, s}^{\tilde{\Gamma}} = \left( -\ch^{\tilde{B}}_3 + s \alpha^2 \tilde{H}^2 \ch^{\tilde{B}}_1 + \tilde{\Gamma} \ch_1^{\tilde{B}} \right) + i \left(\tilde{H} \ch^{\tilde{B}}_2 - \frac{\alpha^2}{2} \tilde{H}^3 \ch^{\tilde{B}}_0 \right).
\]

The goal of this section is to construct a Bridgeland stability condition
$(\AA^{\alpha, \beta}(\tilde{X}), Z_{\alpha, \beta, s}^{\tilde{\Gamma}})$ on $\Db(\tilde{X})$. In order to do so, we need to better understand the relationship between $\Db(X)$ and $\Db(\tilde{X})$. Let $\EE = \OO_{\tilde{X}} \oplus \OO_{\tilde{X}}(E) \oplus \OO_{\tilde{X}}(2E)$, and let $\BB$ be the sheaf of $\OO_X$-algebras given by $f_* \lEnd(\EE)$. We denote the category of finitely generated $\BB$-modules by $\Coh(X, \BB)$ and its bounded derived category by $\Db(X, \BB)$. Toda proves the following theorem:

\begin{thm}[{\cite[Theorem 4.5]{Tod13:extremal_contractions}}]
\label{thm:per_equivalence}
The map $\Phi: Rf_* R\lHom(\EE, - ): \Db(\tilde{X}) \to \Db(X, \BB)$ is an equivalence.
%that restricts to an equivalence between $\Per(\tilde{X}/X)$ and $\Coh{\AA}$.
\end{thm}

\begin{lem}
As sheaves of $\OO_X$-modules there is an isomorphism
\begin{equation}
\label{eq:algebra_A}
\BB \cong \II_Z \oplus \II_P^{\oplus 2} \oplus \OO_X^{\oplus 6},
\end{equation}
where $Z \subset X$ is a zero-dimensional subscheme of length $4$ that is set theoretically supported at $P \in X$.
\end{lem}

\begin{proof}
Note that as $\OO_{\tilde{X}}$-modules
\[
\lEnd(\EE) \cong \OO(-2E) \oplus \OO(-E)^{\oplus 2} \oplus \OO^{\oplus 3} \oplus \OO(E)^{\oplus 2} \oplus \OO(2E).
\]
By applying the derived pushforward $f_*$ to the exact sequence
\[
0 \to \OO \to \OO(E) \to \OO_E(E) = \OO_E(-1) \to 0,
\]
we get $f_* \OO(E) \cong \OO$. We can do the same to the exact sequence
\[
0 \to \OO(E) \to \OO(2E) \to \OO_E(2E) = \OO_E(-2) \to 0,
\]
to get $f_* \OO(2E) \cong \OO$. By using the exact sequence
\[
0 \to \OO(-E) \to \OO \to \OO_E \to 0,
\]
we get $f_* \OO(-E) \cong \II_P$. Lastly, we use the exact sequence
\[
0 \to \OO(-2E) \to \OO(-E) \to \OO_E(-E) = \OO_E(1) \to 0,
\]
to get
\[
0 \to f_* \OO(-2E) \to \II_P \to \OO_P^{\oplus h^0(\OO_E(1))} = \OO_P^{\oplus 3} \to 0.
\]
The claim follows, because $\ch(f_* \OO(-2E)) = (1,0,0,-4)$.
\end{proof}

Instead of working with $\Db(\tilde{X})$ we will construct the stability condition on $\Db(X, \BB)$. First, we have to understand the central charge $Z_{\alpha, \beta, s}^{\tilde{\Gamma}}$ in terms of $\Db(X, \BB)$. Let $j: \Db(X, \BB) \to \Db(X)$ be the forgetful functor. Then we have a Chern character on $\Db(X, \BB)$ given by the composition of $\ch$ and $j$. By abuse of notation we will still call it $\ch$.

\begin{lem}
\label{lem:comparision_Per_A}
For any object $F \in \Db(\tilde{X})$ the equality 
\[
\ch^B(\Phi(F)) = 3 f_* \left( \ch^{\tilde{B}}(F) \cdot \left(1 + \frac{E^2}{6} \right) \right)
\]
holds. In particular, we get $3 Z_{\alpha, \beta, s}^{\tilde{\Gamma}}(F) = Z^{\Gamma}_{\alpha, \beta, s}(\Phi(F))$.
\begin{proof}
Let $Q$ be the universal quotient bundle of $E = \P(\NN_{P/X}) = \P(\OO_P^{\oplus 3}) = \P^2$ and $i \colon E\hookrightarrow \tilde{X}$ the inclusion. By \cite[Lemma 15.4]{Ful98:intersection_theory} we have the following exact columns and row. The exactness of the row is precisely part (i) of that lemma, while the exactness of the third column is part (iii).
% \[
% 0 \to T_{\tilde{X}} \to f^* T_X \to i_* Q \to 0
% \]
\[
\xymatrix{
& 0 & 0 & 0 & \\
0 \ar[r] & i_* \OO_E(E) \ar[r] \ar[u] & i_* \OO_E^{\oplus 3} \ar[r] \ar[u] & i_* Q \ar[u]
\ar[r] & 0
\\
& \OO_{\tilde{X}}(E) \ar[u] & \OO_{\tilde{X}}^{\oplus 3} \ar[u] & f^* T_X \ar[u] & \\
& \OO_{\tilde{X}} \ar[u] & \OO_{\tilde{X}}(-E)^{\oplus 3} \ar[u] & T_{\tilde{X}} \ar[u] & \\
& 0 \ar[u] & 0 \ar[u] & 0 \ar[u] &
}
\]
We used the fact that $\OO_E(E) \cong \OO_E(-1)$. The Grothendieck-Riemann-Roch Theorem implies
\[
\ch(Rf_* R\lHom(\EE, F)) = f_*(\ch (R\lHom(\EE, F)) \cdot \td^{-1}(i_* Q)).
\]
Moreover, we can compute
\begin{align*}
\td^{-1}(i_* Q) &= \td^{-1}(\OO_E^{\oplus 3}) \td(\OO_E(E)) \\
&= \td^{-1}(\OO_{\tilde{X}}^{\oplus 3}) \td(\OO_{\tilde{X}}(-E)^{\oplus 3})
\td(\OO_{\tilde{X}}(E)) \td^{-1}(\OO_{\tilde{X}}) \\
&= \left(1-\frac{E}{2}+\frac{E^2}{12}\right)^3 \left(1 + \frac{E}{2} +
\frac{E^2}{12}\right)
\\
&= 1 - E + \frac{E^2}{3}.
\end{align*}
Therefore, we obtain
\begin{align*}
\ch^B(\Phi(F)) &= f_* \left(\ch^{f^* B} (R\lHom(\EE, F)) \cdot \td^{-1}(i_*
Q)\right)
\\
&= f_*\left( \ch^{f^*B} (F) \cdot \ch(\EE^{\vee}) \cdot \left(1 - E +
\frac{E^2}{3}\right) \right) \\
&= f_* \left( \ch^{f^*B} (F) \cdot \left(3 - 3E + \frac{5 E^2}{2} -
\frac{3E^3}{2} \right) \cdot \left(1 - E + \frac{E^2}{3}\right) \right) \\
&= 3 f_* \left(\ch^{f^*B} (F) \cdot \left(1 - 2E + \frac{13}{6} E^2 -
\frac{5}{3} E^3 \right) \right) \\
&= 3 f_* \left(\ch^{\tilde{B}} (F) \cdot \left(1 + \frac{E^2}{6} \right) \right). \qedhere
\end{align*}
\end{proof}
\end{lem}

We can define the Chern character of a coherent $\BB$-module $F$ via the forgetful functor by declaring $\ch(F) := \ch(j(F))$. Thus, slope stability on $\Coh(X,\BB)$ can be defined with the same slope function as before. Instead of asking for subsheaves to have smaller slope, we ask for $\BB$-submodules to have smaller slope. Since $j$ is faithful, the fact that $\Coh(X)$ is noetherian implies $\Coh(X, \BB)$ to be noetherian. Thus, standard arguments \cite[Proposition 4.10]{MS16:lectures_notes} imply $\mu$ to have Harder-Narasimhan filtrations in $\Coh(X, \BB)$ as well.

The strategy from here on is to construct a double tilt of $\Coh(X, \BB)$ similarly to the case of $\Coh(X)$ in Section \ref{sec:prelim}. Comparing stability via the forgetful functor $j$ will lead to a proof of a Bogomolov-Gieseker type inequality that allows to finish the construction. The following lemmas show that the forgetful functor $j$ not only maps the heart $\Coh(X, \BB)$ to $\Coh(X)$, but does the same to further tilts as well.

\begin{lem}
\label{lem:mu_preserved_forgetful}
Let $F \in \Coh(X, \BB)$ be $\mu$-semistable. Then so is $j(F)$.
\end{lem}

\begin{proof}
Tensor products in this proof will not be derived. Let $G \in \Coh(X, \BB)$ be an arbitrary object, and let $T \into j(G)$ be the maximal torsion subsheaf in $\Coh(X)$. We claim that $T$ is also a $\BB$-module, and $T \into G$ is a morphism in $\Coh(X, \BB)$. Indeed, the action of $\BB$ on $G$ commutes with the action of $\OO_X$, and thus, torsion is mapped to torsion.

Assume $D \into j(F)$ to be a morphism in $\Coh(X)$ where $D$ is $\mu$-semistable in $\Coh(X)$ and $\mu(D) > \mu(F)$. This means $\mu(F) < \infty$ and we can assume both $D$ and $j(F)$ to be torsion free. We get a non-trivial morphism $D \otimes \BB \to F$ in $\Coh(X, \BB)$. Let $S \in \Coh(X, \BB)$ be the maximal torsion subsheaf of $D \otimes \BB$. Then we get a non trivial morphism $(D \otimes \BB)/S \to F$ in $\Coh(X, \BB)$ . The proof will proceed in two steps. Firstly, we will show $\mu((D \otimes \BB)/S) = \mu(D)$. Secondly, we get a contradiction by showing $(D \otimes \BB)/S$ to be $\mu$-semistable in $\Coh(X, \BB)$.

We have an exact sequence $0 \to \BB \to \OO^{\oplus 9} \to \OO_P^{\oplus 2} \oplus \OO_Z \to 0$ in $\Coh(X)$. The long exact $\Tor$ sequence shows $S = \Tor^1(D, \OO_P^2 \oplus \OO_Z)$ which is set-theoretically supported at $P$. We also know that $\Tor^i(D, \BB)$ is supported at $P$ for $i>0$. This shows $\mu((D \otimes \BB)/S) = \mu(D \otimes^L \BB) = \mu(D)$.

This argument also gives us an injective morphism $(D \otimes \BB)/S \into D^{\oplus 9}$. A subsheaf of a semistable sheaf of the same slope has to be again semistable.
\end{proof}

Exactly as in Section \ref{sec:prelim} we can construct a tilt $\Coh^{\beta}(X, \BB)$ in $\Db(X, \BB)$. The previous lemma shows that the functor $j$ maps $\Coh^{\beta}(X, \BB)$ to $\Coh^{\beta}(X)$. As before we will look at the slope function $\nu_{\alpha, \beta}$ on these categories. For any $\beta \in \Q$, the category $\Coh^{\beta}(X)$ is noetherian, and we can use $j$ to show that $\Coh^{\beta}(X, \BB)$ is also noetherian. Therefore, tilt stability is well defined with the same arguments as in \cite{BMT14:stability_threefolds}. Again as in Section \ref{sec:prelim} we want to construct a second tilt $\AA^{\alpha, \beta}(X, \BB)$ in $\Db(X, \BB)$. In order to see that the functor $j$ maps $\AA^{\alpha, \beta}(X, \BB)$ to $\AA^{\alpha, \beta}(X)$ we need to compare $\nu_{\alpha, \beta}$-stability.

\begin{lem}
\label{lem:nu_preserved_forgetful}
Let $F \in \Coh^{\beta}(X, \BB)$ be $\nu_{\alpha, \beta}$-semistable. Then so is $j(F)$.
\begin{proof}
Assume $D \into j(F)$ to be an injective morphism in $\Coh^{\beta}(X)$ such that $\nu_{\alpha, \beta}(D) > \nu_{\alpha, \beta}(F)$, and $D$ is $\nu_{\alpha, \beta}$-semistable in $\Coh^{\beta}(X)$. This implies $\nu_{\alpha, \beta}(F) < \infty$. The strategy of this proof is to construct an object in $\Coh^{\beta}(X, \BB)$ out of $D$ that can be used for a contradiction. Let $\HH_{\beta}^i$ be the cohomology functor with respect to $\Coh^{\beta}(X, \BB)$. By Lemma \ref{lem:mu_preserved_forgetful} the cohomology functor with respect to $\Coh^{\beta}(X)$ is the composition of $\HH_{\beta}^i$ with the forgetful functor, and therefore, we will abuse notation by calling it $\HH_{\beta}^i$, too.

We have two morphisms $D \to D \otimes^L \BB \to j(F)$  in $\Coh^{\beta}(X)$ whose composition is our original map $D \to j(F)$. Indeed, replace both $D$ and $j(F)$ to be complexes of locally free sheaves such that $D \to j(F)$ is a morphism of complexes. Then we have $D \otimes^L \BB = D \otimes \BB$ and the morphisms can be simply constructed as morphisms of complexes. In particular, all these morphisms are non trivial.

The morphism $D \otimes^L \BB \to j(F)$ is in fact a morphism $D \otimes^L \BB \to F$ of $\BB$-module complexes. Thus, we have an induced non trivial morphism $\HH_{\beta}^0(D \otimes^L \BB) \to F$ in $\Coh^{\beta}(X, \BB)$. Let $S \into \HH_{\beta}^0(D \otimes^L \BB)$ be the biggest proper non-trivial subobject in the Harder-Narasimhan filtration of $\HH_{\beta}^0(D \otimes^L \BB)$ with respect to $\nu_{\alpha, \beta}$ in $\Coh^{\beta}(X, \BB)$. In particular, $\HH_{\beta}^0(D \otimes^L \BB)/S$ is $\nu_{\alpha, \beta}$-semistable in $\Coh^{\beta}(X, \BB)$.

As in Lemma \ref{lem:mu_preserved_forgetful} the proof will proceed in two steps. Firstly, we prove $\nu_{\alpha, \beta}(\HH_{\beta}^0(D \otimes^L \BB)/S) = \nu_{\alpha, \beta}(D)$. Secondly, we obtain a contradiction by showing that $\HH_{\beta}^0(D \otimes^L \BB)/S \to F$ is non-trivial.

Recall that we have a short exact sequence 
\[
0 \to \BB \to \OO^{\oplus 9} \to \OO_P^{\oplus 2} \oplus \OO_Z \to 0
\]
of $\OO_X$-modules. The tensor product with $D$ leads to a distinguished triangle
\[
D \otimes^L \BB \to D^{\oplus 9} \to D \otimes^L (\OO_P^{\oplus 2} \oplus \OO_Z).
\]
Part of the long exact sequence with respect to $\Coh^{\beta}(X)$ is given by
\[
0 \to \HH_{\beta}^{-1}(D \otimes^L (\OO_P^{\oplus 2} \oplus \OO_Z)) \to \HH_{\beta}^0(D \otimes^L \BB) \to D^{\oplus 9} \to \HH_{\beta}^{0}(D \otimes^L (\OO_P^{\oplus 2} \oplus \OO_Z)) \to \HH_{\beta}^1(D \otimes^L \BB) \to 0.
\]
Since $\OO_P^{\oplus 2} \oplus \OO_Z$ is set-theoretically supported at $P$, so are $\Tor^i(D,(\OO_P^{\oplus 2} \oplus \OO_Z))$ for every $i$. Thus, $\Tor^i(D,(\OO_P^{\oplus 2} \oplus \OO_Z))$ is an object in $\Coh^{\beta}(X)$ for every $i$, implying that 
\[
\HH_{\beta}^{-i}(D\otimes^L(\OO_P^{\oplus 2} \oplus \OO_Z))=\Tor^i(D,(\OO_P^{\oplus 2} \oplus \OO_Z)).
\]
Therefore, from Lemma \ref{tor:quotient} it follows that $\HH_{\beta}^1(D \otimes^L \BB)$ is also set-theoretically supported at $P$, and the equality $\nu_{\alpha, \beta}(D) = \nu_{\alpha, \beta}(\HH_{\beta}^0(D \otimes^L \BB))$ holds.

We want to show that the inclusion $j(S) \into \HH_{\beta}^0(D \otimes^L \BB)$ in $\Coh^{\beta}(X)$ lifts to an inclusion 
\[
j(S) \into \HH_{\beta}^{-1}(D \otimes^L (\OO_P^{\oplus 2} \oplus \OO_Z)) =: G.
\]
The kernel of the composition $j(S) \into \HH_{\beta}^0(D \otimes^L \BB) \to D^{\oplus 9}$ is a subobject $G' \into G$. Since $G$ is set-theoretically supported on $P$, then by Lemma \ref{tor:quotient} so is $G'$. We have an inclusion $j(S)/G' \into D^{\oplus 9}$, and 
\[
\nu_{\alpha, \beta}(j(S)/G') = \nu_{\alpha, \beta}(j(S)) >\nu_{\alpha,\beta}(\HH_{\beta}^0(D\otimes^L \BB))= \nu_{\alpha, \beta}(D) = \nu_{\alpha, \beta}(D^{\oplus 9}).
\]
This is a contradiction to the stability of $D$ unless $j(S) = G'$, i.e., $j(S) \subset G$. Since $F$ is $\nu_{\alpha, \beta}$-semistable in $\Coh^{\beta}(X, \BB)$, there is no morphism from $S$ to
$F$ in $\Coh^{\beta}(X, \BB)$. Therefore, the non trivial morphism $\HH_{\beta}^0(D \otimes^L \BB) \to F$
induces a non trivial morphism $\HH_{\beta}^0(D \otimes^L \BB)/S \to F$ in $\Coh^{\beta}(X, \BB)$. But since $\nu_{\alpha, \beta}(\HH_{\beta}^0(D \otimes^L \BB)/S) = \nu_{\alpha, \beta}(D)$ that is a contradiction to semistability of $F$.
\end{proof}
\end{lem}

The final step is to show that $(\AA^{\alpha, \beta}(X, \BB), Z_{\alpha, \beta, s}^{\Gamma})$ is a Bridgeland stability condition on $\Db(X, \BB)$. We define a heart $\AA^{\alpha, \beta}(\tilde{X}) = \Phi^{-1}\AA^{\alpha, \beta}(X, \BB) \subset \Db(\tilde{X})$.

\begin{thm}
\label{thm:stability_blow_up}
If the pair $(\AA^{\alpha, \beta}(X), Z_{\alpha, \beta, s}^{\Gamma})$ is a Bridgeland stability condition on $\Db(X)$, then $(\AA^{\alpha, \beta}(\tilde{X}), Z_{\alpha, \beta, s}^{\tilde{\Gamma}})$ is a
Bridgeland stability condition on $\Db(\tilde{X})$. %with support property 
\begin{proof}
By Lemma \ref{lem:nu_preserved_forgetful} the forgetful functor maps $\AA^{\alpha, \beta}(X, \BB)$ to $\AA^{\alpha, \beta}(X)$. Since $Z_{\alpha, \beta, s}^{\Gamma} (\AA^{\alpha, \beta}(X) \backslash \{0\})$ is contained in the union of the upper half-plane and the negative real line, the same holds for $(\AA^{\alpha, \beta}(X, \BB), Z_{\alpha, \beta, s}^{\Gamma})$. In \cite{BMT14:stability_threefolds} it was shown that $\AA^{\alpha, \beta}(X)$ is noetherian for $\alpha$ and $\beta$ rational. Therefore, $\AA^{\alpha, \beta}(X, \BB)$ has to be noetherian, too. Thus, standard arguments \cite[Proposition 4.10]{MS16:lectures_notes} imply $Z_{\alpha, \beta, s}^{\tilde{\Gamma}}$ to have Harder-Narasimhan filtrations in $\AA^{\alpha,\beta}(X, \BB)$ even for irrational $\alpha, \beta$. The fact that $(\AA^{\alpha, \beta}(\tilde{X}), Z_{\alpha, \beta, s}^{\tilde{\Gamma}})$ is a Bridgeland stability condition in $\Db(\tilde{X})$ follows by applying the functor $\Phi^{-1}$ and using Lemma \ref{lem:comparision_Per_A}.
\end{proof}
\end{thm}

\subsection{Stability of skyscraper sheaves}

We finish the article by proving the following proposition about stability of skyscraper sheaves.

\begin{prop}
\begin{enumerate}
\item For any $x \in \tilde{X}$ the skyscraper sheaf $\C(x)$ is contained in $\AA^{\alpha, \beta}(\tilde{X})$.
\item If $x \in \tilde{X} \backslash E$, then $\C(x)$ is stable with respect to $Z_{\alpha, \beta, s}^{\tilde{\Gamma}}$.
\item If $x \in E$, then $\C(x)$ is strictly semistable with respect to $Z_{\alpha, \beta, s}^{\tilde{\Gamma}}$.
\end{enumerate}
\end{prop}

\begin{proof}
\begin{enumerate}
\item We have
\[
\Phi(\C(x)) = Rf_* (\C(x)^{\oplus 3}) = \C(f(x))^{\oplus 3}.
\]
This object is both slope semistable and tilt semistable on $X$ even without the $\BB$-module structure. Since it has slope infinity for both the stabilities, it has to be contained in $\AA^{\alpha, \beta}(X, \BB)$.
\item More strongly, we will show that for $x \in \tilde{X} \backslash E$ the skyscraper sheaf $\C(x)$ is simple, i.e., it has no non-trivial subobjects in $\AA^{\alpha, \beta}(\tilde{X})$. We will do this in three steps. First, we show that $\Phi(\C(x))$ is simple in $\Coh(X, \BB)$, then we do the same for $\Coh^{\beta}(X, \BB)$, and finally for $\AA^{\alpha, \beta}(X, \BB)$.

We need to understand the action of $\BB$ on $\Phi(\C(x)) = \C(f(x))^{\oplus 3}$. This action is completely determined by understanding the action of the restriction $\BB_{|f(x)}$. If $U \subset X$ is an open subset containing $f(x)$, but not $P$, then 
\[
\BB(U) = \End(\EE(f^{-1}(U))) = \Hom(\OO_U^{\oplus 3}, \OO_U^{\oplus 3}).
\]
Therefore, the restriction $\BB_{|f(x)}$ is simply the whole endomorphism algebra of $\C(f(x))^{\oplus 3}$, i.e. the algebra of three times three matrices over $\C$. This action is transitive, and therefore, $\C(x)$ is simple as a $\BB$-module.

Let $E \into \Phi(\C(x))$ be a non-zero subobject in the first tilt $\Coh^{\beta}(X, \BB)$, and let $F$ be the cokernel of this map. By definition of $\Coh^{\beta}(X, \BB)$, we get $E \in \Coh(X, \BB)$. Since $\Phi(\C(x))$ is simple in $\Coh(X, \BB)$, the map $E \to \Phi(\C(x))$ must be surjective as $\BB$-modules. This means $F[-1] \in \Coh(X, \BB)$. However, $F[-1]$ is a submodule of $E$ with the same slope, and thus $E \in \Coh^{\beta}(X, \BB)$ implies $F[-1] \in \Coh^{\beta}(X, \BB)$. This is only possible if $F = 0$, and $E \cong \Phi(\C(x))$,

The previous paragraph with $\Coh(X, \BB)$ replaced by $\AA^{\alpha, \beta}(X, \BB)$ shows that $\C(x)$ has to be simple in $\AA^{\alpha, \beta}(X, \BB)$ as well.
\item If $x \in E$, we will analyze the non-trivial morphism $\OO_E(2E) \to \C(x)$. Since $\OO_{\P^2}(-1)$ and $\OO_{\P^2}(-2)$ have trivial sheaf cohomology, we get
\begin{align*}
\Phi(\OO_E(2E)) &= Rf_* (\OO_E \oplus \OO_E(E) \oplus \OO_E(2E)) \\
&= \C(P) \otimes H^0(\OO_{\P^2}) = \C(P).
\end{align*}
This object is stable in all notions of stability on $X$, with or without $\BB$-module structure. In particular, $\OO_E(2E) \in \AA^{\alpha, \beta}(\tilde{X})$. A straightforward computation shows
\begin{align*}
Z_{\alpha, \beta, s}^{\tilde{\Gamma}}(\OO_E(2E)) &= - \frac{1}{3},\\
Z_{\alpha, \beta, s}^{\tilde{\Gamma}}(\C(x)) &= -1.
\end{align*}
If $\C(x)$ is stable, then the morphism $\OO_E(2E) \to \C(x)$ must be surjective, but that would imply $-1 > -\tfrac{1}{3}$. \qedhere
\end{enumerate}
\end{proof}

\def\cprime{$'$} \def\cprime{$'$}

\appendix
\section{Contracting the section of a Weierstra{\ss} threefold\\
 (By Omprokash Das)}

 \begin{thm}
	Let $X$ be a smooth projective $3$-fold over $\mathbb{C}$ with $K_X\sim 0$. Let $p:X\to S$ be an elliptic fibration in Weierstra{\ss} form over a del Pezzo surface $S$ with canonical section $\sigma: S\to X$. Then there exists a birational contraction $\phi:X\to Y$ contracting exactly the section $\sigma(S)\subset X$ to a point in $Y$. Furthermore, $(Y, \Delta=0)$ has canonical singularities.
	\end{thm}

\begin{proof}
Let $\Theta=\sigma(S)$. From \cite[Corollary 2.2(i)]{Dia16} it follows that $(\Theta-(1+\frac{1}{t}p^*K_S))$ is ample for all $t>0$. Therefore by taking limit as $t\to+\infty$ we see that $(\Theta-p^*K_S)$ is a nef divisor on $X$. By adjunction we have $(K_X+\Theta)|_\Theta\sim K_\Theta$, i.e., $\Theta|_\Theta\sim K_\Theta$ (since $K_X\sim 0$). Note that $-K_\Theta$ is ample, since $\Theta$ is isomorphic to $S$. Then $(\Theta-p^*K_S)^3=\Theta^3=(K_\Theta)^2>0$. Thus we have a nef divisor $(\Theta-p^*K_S)$ on $X$ such that $(\Theta-p^*K_S)^3>0$, hence by \cite[Theorem 2.2.16]{Laz04a} $(\Theta-p^*K_S)$ is a big divisor.
	
Now we will show that the divisor $(\Theta-p^*K_S)$ contracts the divisor $\Theta$ to a point and does not contract anything outside of $\Theta$. To that end note that $(\Theta-p^*K_S)-K_X\sim(\Theta-p^*K_S)$ is nef and big and $X$ is smooth, thus by the base-point free theorem \cite[Theorem 3.3]{KM98} $|m(\Theta-p^*K_S)|$ is a base-point free linear system for $m\gg 0$. In particular, a curve $C\subset X$ is contracted by the morphism associated to the linear system $|m(\Theta-p^*K_S)|$ if and only if $(\Theta-p^*K_S)\cdot C=0$.

Let $C$ be a curve contained in $\Theta\subset X$. Then $(\Theta-p^*K_S)\cdot C=(\Theta|_\Theta)\cdot C-K_\Theta\cdot C=0$ (since $\Theta|_\Theta\sim K_\Theta$). This shows that the divisor $(\Theta-p^*K_S)$ contracts the section $\Theta$ to a point (since it contracts every curve in $\Theta$). Now we will show that it does not contract anything else.
	
Let $\gamma$ be a curve on $X$ contracted by $(\Theta-p^*K_S)$, i.e., $(\Theta-p^*K_S)\cdot\gamma=0$. Then by \cite[Lemma 2.1]{Dia16} there exists a curve $C$ on $S$ such that $\gamma\equiv af+b\sigma_*(C)$, where $f$ is a smooth fiber of $p$. It then follows that $a=0$ and $b\neq 0$, i.e., $\gamma\equiv b\sigma_*(C)$. Note that if $C'\subset X$ is a curve not contained in $\Theta$ then $\Theta\cdot C'\geq 0$; on the other hand $\Theta\cdot\gamma=b(K_\Theta\cdot\sigma_*(C))<0$. Therefore $C'\not\equiv \lambda\gamma$ for any $\lambda\in\mathbb{R}$, and consequiently $(\Theta-p^*K_S)$ does not contract any curve which is not contained in $\Theta$. In other words, $(\Theta-p^*K_S)$ gives a birational divisorial contraction, say $\phi:X\to Y$ such that the exceptional locus of $\phi$ is $\Theta$ and $\phi(\Theta)=\rm{pt}$. Since $X$ and $\Theta$ are both smooth and $\Theta$ is irreducible, $(X, \frac{1}{2}\Theta)$ is a Kawamata log terminal pair (since it is a simple normal crossing pair with coefficients of the boundary divisor in the interval $(0, 1)$). Note that $K_Y=\phi_*K_X\sim 0$. We also see that if $(\Theta-p^*K_S)\cdot C=0$ for a curve $C\subset X$ then $(K_X+\frac{1}{2}\Theta)\cdot C<0$, i.e., $-(K_X+\frac{1}{2}\Theta)$ is $\phi$-nef (it is in fact $\phi$-ample). Then by \cite[Corollary 3.38]{KM98} $(Y, 0)$ has Kawamata log terminal singularities. Now since $K_Y$ is a Cartier divisor ($K_Y\sim 0$), for any exceptional divisor $E$ over $Y$ the discrepancy $a(E, Y)$ is an integer such that $a(E, Y)>-1$, and hence $a(E, Y)\geq 0$. Therefore $(Y, \Delta=0)$ has canonical singularities.     
\end{proof}


\begin{thebibliography}{BBHR09}

\bibitem[BBHR09]{BBH09:fm_nahm_transform}
C.~Bartocci, U.~Bruzzo, and D.~Hern\'andez~Ruip\'erez.
\newblock {\em Fourier-{M}ukai and {N}ahm transforms in geometry and
  mathematical physics}, volume 276 of {\em Progress in Mathematics}.
\newblock Birkh\"auser Boston, Inc., Boston, MA, 2009.

\bibitem[BMS16]{BMS14:abelian_threefolds}
A.~Bayer, E.~Macr{\`{\i}}, and P.~Stellari.
\newblock The space of stability conditions on abelian threefolds, and on some
  {C}alabi-{Y}au threefolds.
\newblock {\em Invent. Math.}, 206(3):869--933, 2016.

\bibitem[BMSZ17]{BMSZ16:stability_fano}
M.~Bernardara, E.~Macr\`\i, B.~Schmidt, and X.~Zhao.
\newblock Bridgeland stability conditions on {F}ano threefolds.
\newblock {\em \'Epijournal Geom. Alg\'ebrique}, 1:Art. 2, 24, 2017.

\bibitem[BMT14]{BMT14:stability_threefolds}
A.~Bayer, E.~Macr{\`{\i}}, and Y.~Toda.
\newblock Bridgeland stability conditions on threefolds {I}:
  {B}ogomolov-{G}ieseker type inequalities.
\newblock {\em J. Algebraic Geom.}, 23(1):117--163, 2014.

\bibitem[Bri07]{Bri07:stability_conditions}
T.~Bridgeland.
\newblock Stability conditions on triangulated categories.
\newblock {\em Ann. of Math. (2)}, 166(2):317--345, 2007.

\bibitem[BvdB03]{BvdB03:functors}
A.~Bondal and M.~van~den Bergh.
\newblock Generators and representability of functors in commutative and
  noncommutative geometry.
\newblock {\em Mosc. Math. J.}, 3(1):1--36, 258, 2003.

\bibitem[CH18]{CH15:nef_cones}
I.~Coskun and J.~Huizenga.
\newblock The nef cone of the moduli space of sheaves and strong {B}ogomolov
  inequalities.
\newblock {\em Israel J. Math.}, 226(1):205--236, 2018.

\bibitem[Dia16]{Dia16:fm_elliptic_CY3}
D.~E. Diaconescu.
\newblock Vertical sheaves and {F}ourier-{M}ukai transform on elliptic
  {C}alabi-{Y}au threefolds.
\newblock {\em Commun. Number Theory Phys.}, 10(3):373--431, 2016.

\bibitem[Ful98]{Ful98:intersection_theory}
W.~Fulton.
\newblock {\em Intersection theory}, volume~2 of {\em Ergebnisse der Mathematik
  und ihrer Grenzgebiete. 3. Folge. A Series of Modern Surveys in Mathematics
  [Results in Mathematics and Related Areas. 3rd Series. A Series of Modern
  Surveys in Mathematics]}.
\newblock Springer-Verlag, Berlin, second edition, 1998.

\bibitem[HRS96]{HRS96:tilting}
D.~Happel, I.~Reiten, and S.~O. Smal{\o}.
\newblock Tilting in abelian categories and quasitilted algebras.
\newblock {\em Mem. Amer. Math. Soc.}, 120(575):viii+ 88, 1996.

\bibitem[Kos18]{Kos17:stability_products}
N.~Koseki.
\newblock Stability conditions on product threefolds of projective spaces and
  {A}belian varieties.
\newblock {\em Bull. Lond. Math. Soc.}, 50(2):229--244, 2018.

\bibitem[Li15]{Li15:conjecture_fano_threefold}
C.~Li.
\newblock Stability conditions on {F}ano threefolds of {P}icard number one,
  2015.
\newblock arXiv:1510.04089.

\bibitem[Mac14a]{Mac14:nested_wall_theorem}
A.~Maciocia.
\newblock Computing the walls associated to {B}ridgeland stability conditions
  on projective surfaces.
\newblock {\em Asian J. Math.}, 18(2):263--279, 2014.

\bibitem[Mac14b]{Mac14:conjecture_p3}
E.~Macr{\`{\i}}.
\newblock A generalized {B}ogomolov-{G}ieseker inequality for the
  three-dimensional projective space.
\newblock {\em Algebra Number Theory}, 8(1):173--190, 2014.

\bibitem[MP16]{MP16:conjecture_abelian_threefoldsII}
A.~Maciocia and D.~Piyaratne.
\newblock Fourier--{M}ukai transforms and {B}ridgeland stability conditions on
  abelian threefolds {II}.
\newblock {\em Internat. J. Math.}, 27(1):1650007, 27, 2016.

\bibitem[MS17]{MS16:lectures_notes}
E.~Macr\`i and B.~Schmidt.
\newblock Lectures on {B}ridgeland stability.
\newblock In {\em Moduli of curves}, volume~21 of {\em Lect. Notes Unione Mat.
  Ital.}, pages 139--211. Springer, Cham, 2017.

\bibitem[Piy17]{Piy17:Fano}
D.~Piyaratne.
\newblock Generalized {B}ogomolov-{G}ieseker type inequalities and {F}ano
  3-folds, 2017.
\newblock arXiv:1705.04011.

\bibitem[Sch14]{Sch14:conjecture_quadric}
B.~Schmidt.
\newblock A generalized {B}ogomolov-{G}ieseker inequality for the smooth
  quadric threefold.
\newblock {\em Bull. Lond. Math. Soc.}, 46(5):915--923, 2014.

\bibitem[Sch17]{Sch17:counterexample}
B.~Schmidt.
\newblock Counterexample to the {G}eneralized {B}ogomolov--{G}ieseker
  {I}nequality for {T}hreefolds.
\newblock {\em Int. Math. Res. Not. IMRN}, (8):2562--2566, 2017.

\bibitem[Tod13]{Tod13:extremal_contractions}
Y.~Toda.
\newblock Stability conditions and extremal contractions.
\newblock {\em Math. Ann.}, 357(2):631--685, 2013.

\end{thebibliography}

\begin{thebibliography}{KM98}

\bibitem[Dia16]{Dia16}
D.~E. Diaconescu, \textsl{ Vertical sheaves and {F}ourier-{M}ukai transform on
  elliptic {C}alabi-{Y}au threefolds},
\newblock Commun. Number Theory Phys. \textbf{ 10}(3), 373--431 (2016).

\bibitem[KM98]{KM98}
J.~Koll{\'a}r and S.~Mori,
\newblock \textsl{ Birational geometry of algebraic varieties}, volume 134 of
  \textsl{ Cambridge Tracts in Mathematics},
\newblock Cambridge University Press, Cambridge, 1998,
\newblock With the collaboration of C. H. Clemens and A. Corti, Translated from
  the 1998 Japanese original.

\bibitem[Laz04]{Laz04a}
R.~Lazarsfeld,
\newblock \textsl{ Positivity in algebraic geometry. {I}}, volume~48 of
  \textsl{ Ergebnisse der Mathematik und ihrer Grenzgebiete. 3. Folge. A Series
  of Modern Surveys in Mathematics [Results in Mathematics and Related Areas.
  3rd Series. A Series of Modern Surveys in Mathematics]},
\newblock Springer-Verlag, Berlin, 2004,
\newblock Classical setting: line bundles and linear series.

\end{thebibliography}
\end{document}